# Cartes auto-organisées pour l'analyse exploratoire de données et la visualisation

**Marie Cottrell, Smaïl Ibbou, Patrick Letrémy, Patrick Rousset.**

SAMOS-MATISSE (UMR 8595), Université Paris 1
90, rue de Tolbiac
75634 PARIS Cedex 13, FRANCE
E-Mail : cottrell@univ-paris1.fr

### Résumé

L'objet de cet article est de présenter comment on peut se servir de l'algorithme de Kohonen pour représenter des données multidimensionnelles, grâce à la propriété d'auto-organisation. On montre comment obtenir ces cartes tant pour des données quantitatives que pour des données qualitatives, ou pour un mélange des deux.

Ce qui est exposé reprend en tout ou en partie différents travaux des membres du SAMOS, en particulier de E. de Bodt, B. Girard, P. Letrémy, S. Ibbou, P. Rousset.

La plupart des exemples présentés ont été étudiés à l'aide des programmes écrits par Patrick Letrémy, en langage IML-SAS, et qui sont disponibles sur le site http://samos.univ-paris1.fr.

### Abstract

This paper shows how to use the Kohonen algorithm to represent multidimensional data, by exploiting the self-organizing property. It is possible to get such maps as well for quantitative variables as for qualitative ones, or for a mixing of both.

The contents of the paper come from various works by SAMOS-MATISSE members, in particular by E. de Bodt, B. Girard, P. Letrémy, S. Ibbou, P. Rousset.

Most of the examples have been studied with the computation routines written by Patrick Letrémy, with the language IML-SAS, which are available on the WEB page http://samos.univ-paris1.fr.

## 1. Introduction

Pour étudier, résumer, représenter des données multidimensionnelles comprenant à la fois des variables quantitatives (à valeurs continues réelles) et qualitatives (discrètes, ordinales ou nominales), les praticiens ont à leur disposition de très nombreuses méthodes performantes, éprouvées et déjà implantées dans la plupart des logiciels statistiques. L'analyse de données consiste à construire des représentations simplifiées de données brutes, pour mettre en évidence les relations, les dominantes, la structure interne du *nuage* des observations. On peut distinguer deux grands groupes de techniques classiques : les *méthodes factorielles* et les *méthodes de classification.*

Les *méthodes factorielles* sont essentiellement linéaires ; elles consistent à chercher des sous-espaces vectoriels, des changements de repères, permettant de réduire les dimensions tout en perdant le moins d'information possible. Les plus connues sont *l'Analyse en Composantes Principales* (ACP) qui permet de projeter des données quantitatives sur les axes les plus significatifs et *l'Analyse des Correspondances* qui permet d'analyser les relations entre les différentes modalités de variables qualitatives croisées (Analyse Factorielle des Correspondances, AFC, pour deux variables, Analyse des Correspondances Multiples, ACM, pour plus de deux variables).

Les *méthodes de classification* sont très nombreuses et diverses. Elles permettent de grouper et de ranger les observations. Les plus utilisées sont la *Classification Hiérarchique* où le nombre de classes n'est pas fixé a priori et la *Méthode des Centres Mobiles* où on cherche à regrouper les données en un certain nombre de classes. La méthode des Centres Mobiles a été étendue aux méthodes dites de *nuées dynamiques.*

On trouvera par exemple dans Lebart et coll. (1995) ou Saporta (1990) une présentation de ces méthodes avec de nombreux exemples. On peut aussi se reporter à des références originelles comme Lebart et coll. (1984), Benzécri (1973) ou Diday et coll. (1983).

Il faut également mentionner les méthodes de *projections révélatrices*, qui cherchent des axes de projection permettant de mettre en évidence les groupes éventuels (ou quelque autre *structure*) : voir Friedman et Tukey (1974) pour les origines, Caussinus et Ruiz-Gazen (2003) pour une présentation résumée récente.

Plus récemment, depuis les années 80, de nouvelles méthodes sont apparues, connues sous le nom de *méthodes neuronales*. Elles proviennent de travaux pluridisciplinaires où se sont retrouvés des biologistes, des physiciens, des informaticiens, des théoriciens du signal, des cogniticiens et, plus récemment encore, des mathématiciens et notamment des statisticiens. Le petit ouvrage de Blayo et Verleysen (1996) dans la collection Que Sais-je est une excellente introduction à ce domaine.

Outre le fait qu'elles sont partiellement issues d'une inspiration biologique ou cognitive, ces méthodes ont rencontré rapidement un certain succès en particulier à cause de leur caractère de « boîte noire », d'outil à tout faire, ayant de très nombreux domaines d'applications. Une fois dépassés un certain excès d'enthousiasme et des difficultés de mise en œuvre, les chercheurs et utilisateurs disposent maintenant d'un arsenal de techniques alternatives, non-linéaires en général et algorithmiques. On pourra consulter par exemple l'ouvrage de Ripley (1996) qui intègre les techniques neuronales parmi les méthodes statistiques. En particulier, les statisticiens commencent à intégrer ces méthodes parmi l'ensemble de leurs outils. Voir à ce sujet les nouveaux modules neuronaux des grands logiciels de Statistique (SAS, SYLAB, STATLAB, S+) et des logiciels de calcul (MATLAB, R, GAUSS, etc.).

Le plus connu des modèles neuronaux reste sans conteste le modèle du *Perceptron Multicouches* (voir par exemple Rumelhart et McClelland., 1986, Le Cun, 1986) dont nous ne parlerons pas ici. On peut consulter les chapitres correspondants des ouvrages classiques comme par exemple Hertz et coll. (1991), ou Haykin (1994). On trouvera dans Rynkiewicz et coll. (2001), une présentation résumée et clairement statistique du Perceptron Multicouches.

Nous présentons ici des méthodes de classification, de représentation, d'analyse des relations, toutes construites à partir du célèbre algorithme de Kohonen : voir par exemple Blayo et Demartines (1991, 1992), Varfis et Versino (1992), Kaski (1997), ou une synthèse dans Cottrell et Rousset (1997). De nombreux exemples sont également présentés dans Deboeck et Kohonen (1998), ou dans l'ouvrage collectif *Kohonen Maps*, publié par Errki Oja et Samuel Kaski en 1999, à l'occasion de la conférence WSOM'99, (Oja et Kaski, 1999).

On trouvera sur le site WEB du SAMOS-MATISSE[1] de nombreuses publications et prépublications comportant des applications de ces méthodes dans des domaines très variés (économie, finance, géographie, gestion,…).

Le présent article est structuré comme suit. Dans le paragraphe 2, nous introduisons les notions sous-jacentes de quantification vectorielle et de classification. Dans le paragraphe 3, nous définissons rapidement l'algorithme de Kohonen qui sert de base à tout le reste. Nous montrons dans le paragraphe 4 comment l'utiliser pour réaliser une classification des observations, basée sur des variables quantitatives et admettant une représentation analogue à une ACP, et comment étudier les classes obtenues. Dans le paragraphe 5, est abordé le problème du traitement des données manquantes. Les paragraphes 6 et 7 sont consacrés à des variantes (algorithme Batch, choix des initialisations). Dans le paragraphe 8, nous présentons une méthode de visualisation des distributions des variables qualitatives sur les classes obtenues.

Les paragraphes suivants traitent des variables qualitatives. Dans les paragraphes 9 à 13, nous définissons plusieurs méthodes analogues à l'Analyse des

---

[1] Adresse : http://samos.univ-paris1.fr

Correspondances double (KORRESP) et multiple (KACM, KACM1, KACM2, KDISJ). Le paragraphe 14 est une conclusion (provisoire bien entendu).

Les exemples présentés proviennent de plusieurs études réelles dont les références détaillées sont citées dans l'Annexe E. L'exemple I concerne des profils horaires de consommation électrique en Pologne. L'exemple II est extrait d'une étude, menée en collaboration avec l'équipe PARIS de l'université Paris 1, sur les communes de la vallée du Rhône. L'exemple III étudie les relations entre la nature des monuments en France et leurs propriétaires. L'exemple IV traite de données issues de l'enquête INSEE *Emploi du temps*, concernant des travailleurs à temps partiel.

## 2. Quantification vectorielle et classification

L'algorithme SOM (Self-Organizing Map), appelé également algorithme de Kohonen, est un algorithme original de classement qui a été défini par Teuvo Kohonen, dans les années 80, à partir de motivations neuromimétiques (cf. Kohonen, 1984, 1995). Dans le contexte d'analyse des données qui nous intéresse ici, il est plus éclairant de présenter l'algorithme de Kohonen comme une généralisation de la version stochastique d'un algorithme extrêmement populaire, qui est l'algorithme des Centres Mobiles ou algorithme de Forgy (Forgy, 1965).

Dans tout ce qui suit, on considère que l'ensemble des données est fini et stocké dans une table ayant $N$ lignes représentant les individus (ou observations) et $p$ colonnes correspondant aux variables décrivant les individus.

Chacune des $N$ observations ($x_1$, $x_2$, ..., $x_N$) est donc décrite par $p$ variables quantitatives. La matrice ($N \times p$) ainsi formée est appelée *matrice des données*. On ajoutera ensuite éventuellement des variables qualitatives.

Les buts des méthodes de classification et de quantification vectorielle sont voisins. Dans le premier cas, il s'agit de regrouper et de ranger les observations en des classes qui soient à la fois les plus homogènes possibles et les plus distinctes les unes des autres possibles. Dans le second cas, on cherche à représenter un grand nombre de données multidimensionnelles au moyen d'un nombre beaucoup plus petit de vecteurs codes de même dimension que les données, de façon à perdre le moins possible d'information.

En fait, les deux démarches sont duales. Si l'on part d'une classification donnée, chaque classe est naturellement représentée par le centre de gravité (centre, centroïde) des observations qu'elle contient, qui devient donc le vecteur code (ou représentant) de cette classe. Réciproquement, un ensemble de vecteurs codes définit un ensemble de classes, formées des observations les plus proches de chaque vecteur code au sens d'une distance donnée. Dans ce cas, le vecteur code n'est pas nécessairement égal au barycentre de sa classe.

**Qualité d'une quantification vectorielle**

Si on note $n$ le nombre de classes, $C_1$, $C_2$,..., $C_n$, les vecteurs codes, et $\Gamma_1$, $\Gamma_2$,..., $\Gamma_n$, les classes correspondantes, le but d'une quantification vectorielle est de déterminer les vecteurs codes qui minimisent l'erreur quadratique, appelée *distorsion* dans ce contexte, définie par

$$V(x_1, x_2, \cdots, x_N, C_1, C_2, \cdots, C_n) = \sum_{i=1}^{n} \sum_{x_j \in \Gamma_i} \left\| x_j - C_i \right\|^2.$$

Cette distorsion est analogue à la somme des carrés intra-classes, mais les vecteurs codes ne sont pas nécessairement les barycentres des classes.

**Qualité d'une classification**

Pour réaliser une bonne classification à nombre de classes fixé $n$, on cherche à minimiser la somme des carrés intra-classes, qui est identique à la distorsion lorsque les vecteurs codes $C_1$, $C_2$,..., $C_n$, sont exactement les centres de gravité des classes. Dans ce cas, on minimise donc

$$SC_{\text{intra}}(\Gamma_1, \Gamma_2, \cdots \Gamma_N) = \sum_{i=1}^{n} \sum_{x_j \in \Gamma_i} \left\| x_j - \overline{\Gamma_i} \right\|^2,$$

où les $\overline{\Gamma_i}$ sont les centres de gravité des classes $\Gamma_i$.

**Algorithme de Forgy**

Une des méthodes de classification les plus utilisées est la méthode de Forgy. Rappelons la brièvement. On se donne a priori le nombre $n$ de classes. On initialise les vecteurs codes $C_1$, $C_2$,..., $C_n$ de façon quelconque et on en déduit les classes $\Gamma_1$, $\Gamma_2$,...,$\Gamma_n$, par la méthode du plus proche voisin. On affecte donc à la classe $\Gamma_i$ l'ensemble des observations plus proches de $C_i$ que de tous les autres vecteurs codes $C_j$, avec $j \neq i$. On recalcule ensuite tous les vecteurs codes en les prenant égaux aux centres de gravité des classes, et ainsi de suite. Il est bien connu que cet algorithme assure la décroissance de la somme des carrés intra-classes, et qu'il converge donc vers un de ses minima (local).

Lorsque les données sont très nombreuses, ou lorsque qu'on veut travailler « en ligne » parce qu'on ne dispose pas de toutes les données à la fois, on utilise souvent une version stochastique de cet algorithme. Pour introduire l'algorithme de Kohonen, on présente ici la version stochastique de l'algorithme de Forgy la plus simple. On l'appelle d'ailleurs Algorithme d'Apprentissage Compétitif Simple (abrégé en SCL).

**Algorithme d'apprentissage compétitif simple**

On part également de vecteurs codes aléatoires. A chaque itération *t*, on tire au hasard une observation *x(t+1)*, on détermine l'indice $i_0$ du vecteur code le plus proche (on dit qu'il est gagnant, on parle de compétition), en posant :

$$i_0(C_1(t), C_2(t), \cdots, C_n(t), x(t+1)) = Arg \min_i \|x(t+1) - C_i(t)\|.$$

Ensuite on modifie le vecteur code $C_{i_0}(t)$ et lui seul en posant

$$\begin{cases} C_{i_0}(t+1) = C_{i_0}(t) + \varepsilon(t)(x(t+1) - C_{i_0}(t)) \\ C_i(t+1) = C_i(t), \text{ pour } i \neq i_0, \end{cases}$$

où le paramètre $\varepsilon(t)$ est positif, constant ou décroissant vers 0.

Cela définit un algorithme stochastique à pas $\varepsilon(t)$. La suite des $\varepsilon(t)$ est le gain du processus. Si la suite des gains vérifie les conditions de Robbins-Monro

$$\sum_t \varepsilon(t) = +\infty \text{ et } \sum_t \varepsilon^2(t) < +\infty,$$

on peut montrer que cet algorithme converge vers un minimum (local) de la somme des carrés intra-classes.

Cet algorithme stochastique SCL est la version stochastique de l'algorithme de Forgy. Il a été défini et étudié par Mac Queen (1967) ou par Hertz et coll. (1991). Réciproquement, l'algorithme de Forgy est l'algorithme en moyenne associé à l'algorithme SCL.

Pour ces deux algorithmes, le minimum trouvé dépend de l'initialisation.

## 3. Algorithme de Kohonen

On peut maintenant définir l'algorithme de Kohonen (Kohonen, 1984, 1993, 1995) comme une simple généralisation de l'algorithme SCL, en y rajoutant une notion de voisinage entre les vecteurs codes (ou, ce qui est équivalent, entre les classes).

L'algorithme de Kohonen est alors vu comme un algorithme stochastique de classification, qui regroupe les observations en classes, tout en *respectant la topologie de l'espace des observations*.

Cela veut dire que, si l'on définit a priori une notion de voisinage entre classes, des observations voisines dans l'espace des observations (de dimension *p*) appartiennent (après classement) à la même classe *ou à des classes voisines*.

Les voisinages entre classes peuvent être choisis de manière variée, mais on suppose en général que les classes sont disposées sur une *grille* rectangulaire qui définit naturellement les voisins de chaque classe. On peut aussi considérer une topologie unidimensionnelle dite en *ficelle*, ou éventuellement un tore ou un cylindre. On peut aussi considérer des voisinages hexagonaux qui présentent l'avantage de rendre égales les distances apparentes dans toutes les directions.

Par exemple, sur une grille, on pourra prendre des voisinages de rayon 3 (49 voisins), de rayon 2 (25 voisins) ou de rayon 1 (9 voisins).

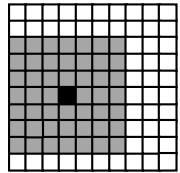 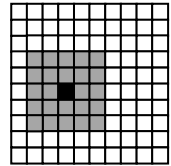 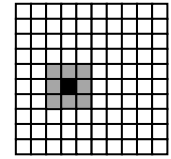

Voisinage de 49   Voisinage de 25   Voisinage de 9

Pour une ficelle, les mêmes rayons donnent 7, 5 et 3 voisins.

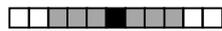 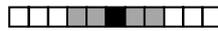 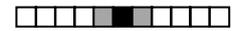

Voisinage de 7   Voisinage de 5   Voisinage de 3

Les classes situées sur les bords ont éventuellement moins de voisines.

**Principe de l'algorithme de Kohonen**

L'algorithme de classement est itératif. L'initialisation consiste à associer à chaque classe un *vecteur code* (ou *représentant*) de *p* dimensions choisi de manière aléatoire. Ensuite, à chaque étape, on choisit une observation au hasard, on la compare à tous les vecteurs codes, et on détermine la classe gagnante, c'est-à-dire celle dont le vecteur code est le plus proche au sens d'une distance donnée a priori. On rapproche alors de l'observation les codes de la classe gagnante *et des classes voisines*.

On voit bien qu'il s'agit d'une généralisation immédiate de l'algorithme SCL. On dit quelquefois que l'algorithme SCL est un algorithme de Kohonen à 0 voisin (en fait il a été défini au moins 20 ans auparavant : voir Mac Queen, 1967).

**Notations et définitions**

On se donne un réseau de Kohonen formé de $n$ unités rangées suivant une certaine topologie (en général une grille, ou bien une ficelle). Ces $n$ unités sont munies d'un système de voisinages homogènes dans l'espace. Pour chaque unité $i$ du réseau, on définit un voisinage de rayon $r$ noté $V_r(i)$ et formé de l'ensemble des unités situées sur le réseau à une distance inférieure ou égale à $r$.

Chaque unité $i$ est représentée dans l'espace $R^p$ par un vecteur $C_i$ (appelé vecteur poids par Kohonen) que nous désignerons par *vecteur code* de l'unité $i$ (ou de la classe $i$). L'état du réseau au temps $t$ est donné par $C(t) = (C_1(t), C_2(t), ..., C_n(t))$.

Pour un état donné $C$ et une observation donnée $x$, l'unité (ou classe) *gagnante* $i_0(C,x)$ est l'unité dont le vecteur code $C_{i_0(C,x)}$ est le plus proche de l'observation $x$ au sens d'une certaine distance. On a donc

$$i_0(C, x) = Arg \min_i \|x - C_i\|.$$

Si l'on pose pour toute unité $i$, $\Gamma_i = \{x \in \{x_1, x_2, ..., x_N\} / i_0(C,x) = i\}$, on dira que $\Gamma_i$ est la classe de numéro $i$ et l'ensemble des classes $(\Gamma_1, \Gamma_2, ..., \Gamma_n)$ forme une partition de l'ensemble des observations, appelée *Partition de Voronoï*. Chaque classe est représentée par le *vecteur code* correspondant. Chaque observation est représentée par le vecteur code le plus proche, exactement comme dans la méthode du plus proche voisin.

Alors, pour un état $C$ donné, le réseau définit une application $\Phi_C$ qui à chaque observation $x$ associe l'unité gagnante correspondante, c'est-à-dire le numéro de sa classe. Après convergence de l'algorithme de Kohonen, l'application $\Phi_C$ *respecte la topologie* de l'espace des entrées, en ce sens que des observations voisines dans l'espace $R^p$ se retrouvent associées à des unités voisines ou à la même unité.

L'algorithme de construction des vecteurs codes est défini de manière itérative comme suit :

- *Au temps 0*, les $n$ vecteurs codes sont initialisés de manière aléatoire (on peut par exemple tirer au hasard $n$ observations).

- *Au temps $t$*, l'état du réseau est $C(t)$ et on présente une observation $x(t+1)$ suivant une loi de probabilité $P$, on a alors :

$$\begin{cases} i_0(C(t), x(t+1)) = Arg \min\{\|x(t+1) - C_i(t)\|, 1 \le i \le n\} \\ C_i(t+1) = C_i(t) - \varepsilon(t)(C_i(t) - x(t+1)), & \forall i \in V_{r(t)}(i_0) \\ C_i(t+1) = C_i(t) & \forall i \notin V_{r(t)}(i_0) \end{cases}$$

où $\varepsilon(t)$ est le *paramètre d'adaptation* ou *de gain,* et où $r(t)$ est le rayon des voisinages au temps $t$.

*Les paramètres importants sont*
- la dimension $p$ de l'espace des entrées,
- la topologie du réseau (grille, ficelle, cylindre, tore, etc.),
- le paramètre d'adaptation, positif, compris entre 0 et 1, constant ou décroissant,
- le rayon des voisinages, en général décroissant,
- la loi de probabilité $P$ de tirage des observations.

On peut écrire cet algorithme sous une forme « classique » d'algorithme stochastique
$$C(t+1) = C(t) + \varepsilon(t)H\bigl(x(t+1), C(t)\bigr).$$

Cependant l'étude mathématique de cet algorithme si facile à définir et à écrire n'est pas simple. Un des points qui rendent difficile l'étude théorique dans le cas général de l'algorithme de Kohonen, est qu'il n'existe pas de potentiel (ou énergie) associé lorsque les entrées sont tirées suivant une loi de probabilité $P$ continue(cf. Erwin et coll., 1992). Mais au contraire, dans le cas qui nous intéresse ici, c'est-à-dire quand l'espace des entrées est fini, Ritter et Schulten (1992) ont montré que l'algorithme de Kohonen pour un rayon constant $r$ de voisinage, est alors un algorithme du gradient stochastique qui minimise le potentiel :

$$V(x_1, x_2, \cdots, x_N, C_1, C_2, \ldots, C_N) = \sum_{i=1}^{n} \sum_{k \in V_r(i)} \sum_{x_j \in \Gamma_k} \left\| x_j - C_i \right\|^2.$$

Ce potentiel généralise la distorsion quadratique ou la variance intra-classes. Ici on calcule la somme des carrés des distances de chaque observation non seulement à son vecteur code mais aussi aux vecteurs codes des classes voisines. On l'appellera dans la suite *Distorsion étendue*.

L'étude mathématique de ce potentiel n'est pas simple non plus, car il n'est pas partout différentiable. En outre, il n'est pas possible d'éviter les minima locaux.

En fait l'étude de la convergence de cet algorithme est donc incomplète et pose des problèmes mathématiques difficiles : voir par exemple Cottrell et Fort (1987), Ritter et Shulten (1986, 1988), Fort et Pagès (1995, 1996, 1999), Benaïm et coll. (1998), Sadeghi (1997, 1998), Cottrell, Fort et Pagès (1995, 1997, 1998). Pour l'instant, l'essentiel des résultats correspond à la dimension 1 (topologie en ficelle et observations de dimension $p = 1$). Des résultats d'organisation et de convergence en loi sont disponibles quand le paramètre d'adaptation est constant. Pour les résultats de convergence presque sûre après réorganisation, il faut que la suite des paramètres $\varepsilon(t)$ vérifie des conditions de Robbins-Monro, classiques pour les algorithmes stochastiques, qui s'énoncent :

$$\sum_t \varepsilon_t = +\infty \quad \text{et} \quad \sum_t \varepsilon_t^2 < +\infty.$$

En d'autres termes, le paramètre doit être « petit », mais « pas trop ».

On trouvera dans Cottrell, Fort et Pagès (1998) un exposé des principaux résultats théoriques disponibles ainsi que des difficultés non résolues.

Un inconvénient de l'algorithme de base est que le nombre de classes doit être fixé a priori. Pour pallier cet inconvénient, on peut, après la classification de Kohonen, pratiquer une classification de type hiérarchique sur les codes des classes de Kohonen, de manière à les regrouper en classes moins nombreuses. Nous allons étudier dans le prochain paragraphe comment obtenir des représentations faciles à interpréter visuellement.

## 4. Classes et super-classes, étude des classes. KACP

Après convergence de l'algorithme, on a vu que les $N$ observations sont classifiées en $n$ classes selon la méthode du plus proche voisin, relativement à la distance choisie dans $R^p$.

### 4.1 Représentations graphiques

On peut alors construire une représentation graphique selon la topologie du réseau. Dans chaque unité du réseau (grille, ficelle, etc.), on dessine les observations associées (superposées) ou on fait la liste de ces observations. Grâce à la propriété de conservation de la topologie, la représentation respecte les relations de voisinage. On obtient ainsi une carte de Kohonen où les grandes caractéristiques du nuage des données sont visibles après un petit nombre d'itérations (de l'ordre de $5N$ ou $6N$, en général). Cette carte est un instrument d'analyse des données qui fournit des informations analogues à celle que donne une Analyse en Composantes Principales (ACP). Bien sûr, il n'y a pas de projection à proprement parler, la carte est grossière. Mais elle est unique, ce qui évite de devoir combiner les différentes projections planes qu'on obtient en ACP. La continuité entre une classe et ses voisines permet de bien comprendre l'évolution le long de la grille, et est facile à interpréter. C'est pour ces raisons que l'ensemble des techniques de classification et représentation fournies par l'algorithme de Kohonen appliqué à une matrice de données est désigné par le sigle KACP.

De nombreux exemples de ces représentations sont disponibles dans la bibliographie ou sur le site WEB du SAMOS-MATISSE : carte des pays repérés par des variables socio-économiques, courbes de consommation électrique, communes d'Ile-de-france, avec des variables immobilières, profils de consommateurs canadiens, segmentation du marché du travail en France, démographie et composition sociale dans la vallée du Rhône, rôle du leasing pour les entreprises de Belgique, classification des ménages en France, etc. Voir Blayo et Demartines (1991), Cottrell, Girard, Girard, Muller et Rousset (1995), Gaubert, Ibbou et Tutin (1996), Cottrell, Girard et Rousset (1998), Cottrell et Gaubert (2000), Gaubert et Cottrell (1999, 1999), Cottrell, Gaubert, Letrémy et Rousset (1999), Ponthieux et Cottrell (2001), etc.

On montre par exemple dans le figure 1, la représentation du contenu de 100 classes, pour une topologie en cylindre 10 par 10, et des données (transformées) de consommations électriques demi-horaires (Exemple I). Dans ce cas, les observations, qui sont des vecteurs de $p = 48$ points, sont dessinées comme des courbes. Dans d'autres cas, on choisira des représentations en histogrammes ou autres.

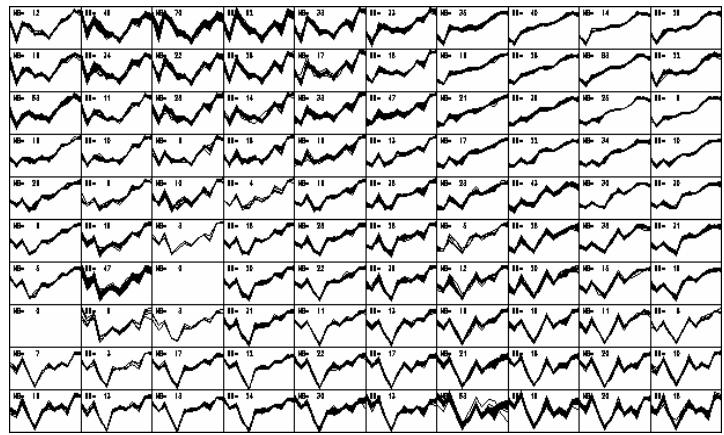

Fig. 1 : *Consommations électriques demi-horaires : contenu des classes.*

Cette représentation permet de contrôler d'un coup d'œil l'homogénéité des classes, les données aberrantes (extrêmes ou erronées) se remarquent tout de suite, etc. Il est également intéressant de représenter les vecteurs codes sur le réseau comme on le voit à la figure 2.

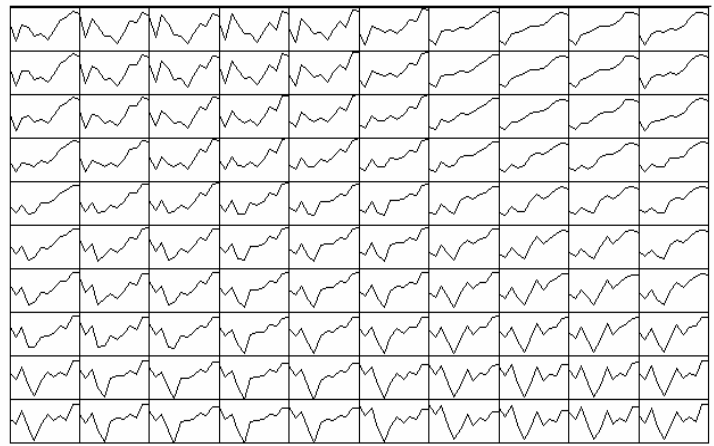

Fig. 2 : *Consommations électriques demi-horaires : vecteurs codes.*

Sur cette représentation, on peut voir immédiatement comment s'est faite l'organisation, les vecteurs codes se déforment de proche en proche sur la carte (rappelons qu'ici il s'agit d'une topologie en cylindre, ce qui explique que les vecteurs codes du bord gauche et du bord droit se ressemblent.

**4.2 Distances entre les classes. Super-classes**

On peut mettre en évidence visuellement les *distance*s entre les classes, qui sont artificiellement équidistantes dans les représentations ci-dessus. Pour cela, en suivant la méthode proposée par de Bodt et Cottrell (1996), on dessine dans chaque cellule un octogone. Dans chacune des 8 directions principales, son sommet est d'autant plus proche du bord que la distance au voisin dans cette direction est petite. Cela permet de faire apparaître les groupes de classes proches et donne une idée de la discrimination entre classes.

Comme le choix du nombre $n$ de classes est arbitraire (et souvent élevé puisqu'on choisit couramment des grilles 8 par 8, 10 par 10 ou plus), on peut réduire le nombre de classes, en les regroupant au moyen d'une *classification hiérarchique* classique sur les $n$ vecteurs codes. On peut alors *colorier* ou *hachurer* les groupes de classes (appelés *super-classes*) pour les rendre visibles. On constate toujours que les super-classes ne regroupent que des classes contiguës, ce qui s'explique par la propriété de respect de la topologie de l'algorithme de Kohonen. D'ailleurs, le non-respect de cette propriété serait un signe de manque de convergence de l'algorithme ou d'une structure particulièrement « repliée » du nuage des données.

Les figures 3 et 4 montrent les deux classements emboîtés (100 classes de Kohonen et 10 super-classes) et les distances inter-classes.

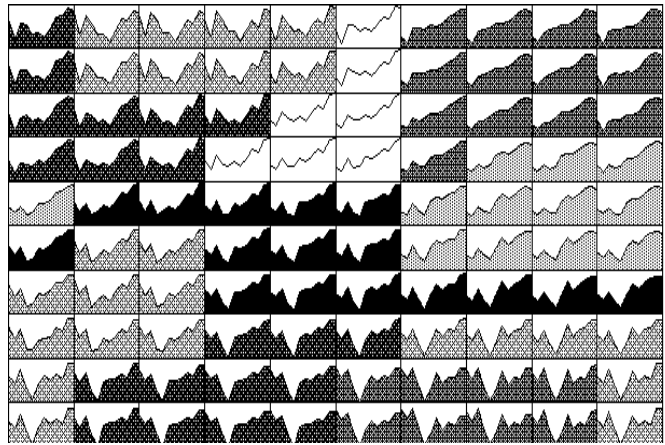

Fig. 3 : *Consommations électriques demi-horaires :
Classes et super-classes «coloriées».*

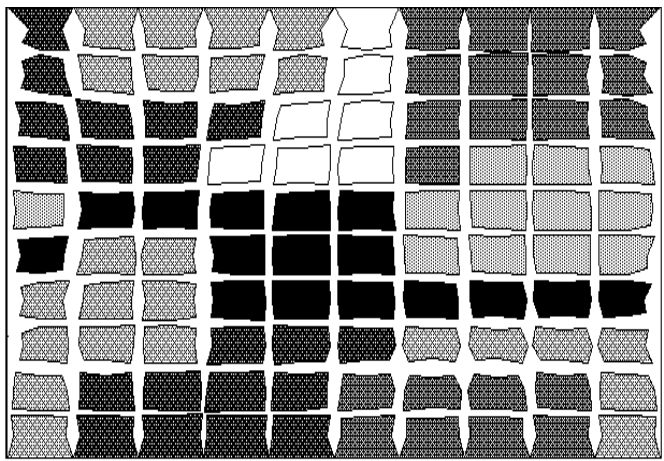

Fig. 4 : *Consommations électriques demi-horaires : Distances inter-classes. Plus il y a de vide entre deux octogones, plus les classes correspondantes sont distantes.*

L'avantage de cette double classification est la possibilité d'analyser les données à un niveau plus grossier, où seules les caractéristiques principales peuvent apparaître, en facilitant l'interprétation.

### 4.3 Typologie

On peut alors établir une typologie des individus en décrivant chacune des super-classes au moyen des statistiques classiques. On peut également, selon le problème, définir un modèle (de régression, d'auto-régression, une analyse factorielle, etc.), adapté à chacune des super-classes.

## 5. Données manquantes

Un des intérêts de l'algorithme de Kohonen utilisé pour la classification est qu'il supporte parfaitement la présence de données manquantes. Voir à ce sujet la thèse de Smaïl Ibbou (1998). Le principe est très simple.

Lorsqu'on présente un vecteur de données incomplet $x$, on calcule le vecteur code gagnant en posant

$$i_0(C, x) = Arg \min_i \|x - C_i\|$$

où la distance $\|x - C_i\|$ est calculée sur les composantes présentes dans le vecteur $x$.

On peut utiliser les vecteurs avec données manquantes de deux façons. Si l'on souhaite les utiliser au moment de la construction des vecteurs-codes, à chaque étape, une fois déterminé le numéro de l'unité gagnante, la mise à jour des vecteurs codes (le gagnant et ses voisins) ne porte que sur les composantes présentes dans le vecteur. On peut aussi se contenter de classer, après construction de la carte, les vecteurs incomplets en les affectant dans la classe dont le vecteur code est le plus proche, au sens de la distance restreinte aux composantes présentes.

Cela donne d'excellents résultats, dans la mesure bien sûr où une variable n'est pas complètement absente ou presque, et aussi dans la mesure où les variables sont corrélées, ce qui est le cas dans la plupart des jeux de données réelles.

On pourra voir de nombreux exemples dans la thèse de Smaïl Ibbou (1998) et aussi dans Gaubert, Ibbou et Tutin (1996).

## 6. Algorithme Batch associé

On a vu que les deux algorithmes SCL et SOM sont des algorithmes stochastiques adaptatifs, ce qui signifie qu'à chaque présentation d'une nouvelle observation, l'algorithme met à jour les valeurs des vecteurs-codes des classes. Mais de même que l'algorithme de Forgy est la version déterministe de l'algorithme SCL, on peut définir un algorithme déterministe équivalent de l'algorithme de Kohonen, qui utilise toutes les données à chaque étape. Il est défini comme un algorithme en moyenne sur l'ensemble des données à chaque étape. On parle d'algorithme de Kohonen Batch, abrégé en KBATCH (Kohonen, 1999).

On sait que l'étude de la fonction *distorsion étendue* ne suffit pas pour montrer la convergence de l'algorithme de Kohonen. Cependant on peut remarquer que, si l'algorithme converge, il converge vers l'un des points d'équilibre de l'équation différentielle ordinaire (ODE) associée à l'algorithme. Cette équation s'écrit avec les notations précédentes (en supposant que la taille du voisinage reste fixe, égale à *r*) :

$$\frac{dC_i(u)}{du} = \sum_{j \in V_r(i)} \sum_{x_l \in \Gamma_j} (C_i(u) - x_l), \ \forall i = 1, \cdots, n$$

Les points d'équilibre de cette équation vérifient alors

$$C_i^* = \frac{\sum_{j \in V_r(i)} \sum_{x_l \in \Gamma_j} x_l}{Card(\bigcup_{j \in V_r(i)} \Gamma_j)}$$

Dans le cas de l'algorithme SCL (où le rayon *r* est égal à 1), ceci exprime que chaque $C_i^*$ est le centre de gravité de sa classe. Dans le cas de SOM, $C_i^*$ est le centre de gravité de l'union de sa classe $\Gamma_i$ et des classes voisines.

On peut alors chercher à calculer ces solutions par un calcul itératif déterministe.

On en déduit immédiatement la définition de l'algorithme KBATCH de calcul des $C_i{}^*$. Ces points sont obtenus comme limite de la suite définie par :

$$C_i^{k+1} = \frac{\sum\limits_{j \in V_r(i)} \sum\limits_{x_l \in \Gamma_j^k} x_l}{Card(\bigcup\limits_{j \in V_r(i)} \Gamma_j^k)},$$

où l'exposant $k$ indique l'étape d'itération et la valeur courante des vecteurs codes et des classes.

Quand on ne prend pas en compte les voisinages, on retrouve exactement l'algorithme de Forgy (les centres mobiles).

On sait (Fort et coll., 2001) que cet algorithme KBATCH est en fait un algorithme "quasi-Newtonien" (approximation de l'algorithme du gradient du second ordre qui minimise la distorsion étendue $V$). C'est un algorithme "quasi-Newtonien" parce qu'on utilise seulement la diagonale de la matrice Hessienne. Malheureusement, il existe de nombreux ensembles disjoints où cette fonction $V$ est différentiable et il y a un minimum local à l'intérieur de chacun d'eux.

On sait aussi que les algorithmes de Newton ne sont pas nécessairement des algorithmes de descente, ce qui entraîne que la fonction distorsion étendue peut croître pour certaines itérations. Cette propriété, qui est un inconvénient, peut néanmoins s'avérer utile, car elle permet ensuite éventuellement la convergence vers un meilleur minimum.

Ces algorithmes déterministes sont, on le sait, extrêmement sensibles aux choix des conditions initiales. On constate que cet algorithme KBATCH, qui en principe minimise la même fonction que l'algorithme SOM, peut conduire à des états d'équilibre mal organisés, alors que l'algorithme stochastique SOM organise parfaitement. On trouvera de nombreux exemples dans Fort et coll. (2002).

Cela conduit à conclure que, si l'algorithme KBATCH peut être séduisant par la rapidité de sa convergence et la simplicité de sa mise en oeuvre, le choix des conditions initiales s'avère crucial du point de vue de la qualité de l'organisation, et qu'il est donc primordial de partir de conditions initiales choisies avec soin. C'est ce que nous allons illustrer dans le paragraphe suivant.

## 7. Diverses initialisations

On peut choisir en théorie les vecteurs codes initiaux de manière quelconque. Le bon sens demande qu'on les choisisse au moins dans l'enveloppe convexe des observations, puisque les vecteurs codes doivent être après convergence des représentants des observations. Cette façon de faire donne de bons résultats dans la plupart des cas, lorsqu'il n'y a pas trop de données excentrées par rapport au reste du nuage (méthode I). On peut aussi choisir (si le nombre de données est grand) *n* observations elles-mêmes prises au hasard comme vecteurs codes (méthode II).

Enfin, et c'est ce qui recommandé, même si cela demande un peu plus de calculs, on peut choisir les vecteurs codes formant un maillage régulier de l'enveloppe convexe des projections des observations sur le premier plan principal (méthode III). Cela accélère la convergence et fournit une meilleure organisation.

Les figures 5, 6, 7 correspondent à une classification des communes de la vallée du Rhône[2] (Exemple II). Il y a 1783 communes. Pour chacune on connaît les chiffres de population totale mesurée par 7 recensements (1936, 1954, 1962, 1968, 1975, 1982, 1990). On pratique d'abord une classification en 64 classes, qui sont ensuite regroupées en 5 super classes. La figure 5 correspond à l'initialisation I. En haut, elle montre les super classes et les 64 vecteurs codes pour l'algorithme SOM et pour KBATCH. En bas, est représenté le contenu des 5 super classes.

---

[2] Les données nous ont été fournies par les collègues de l'équipe PARIS de l'Université Paris 1.

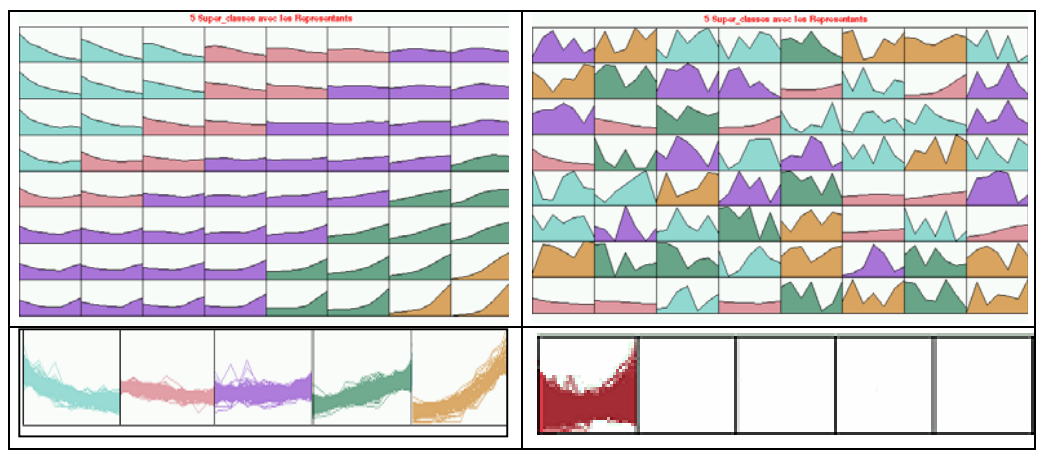

Fig. 5 : *Les communes, les classes et les super-classes,
initialisation I (dans l'enveloppe convexe),
à gauche, algorithme SOM, à droite, algorithme KBATCH*

Pour l'algorithme SOM, le classement est parfait, des communes en forte décroissance dans le coin en haut à gauche, jusqu'aux communes en forte croissance en bas à droite. L'algorithme KBATCH n'organise pas du tout, les super-classes regroupent des classes éloignées, les classes voisines ne se ressemblent pas. De plus toutes les communes se retrouvent dans la même super-classe.

Les figures 6 et 7 correspondent aux initialisations II et III.

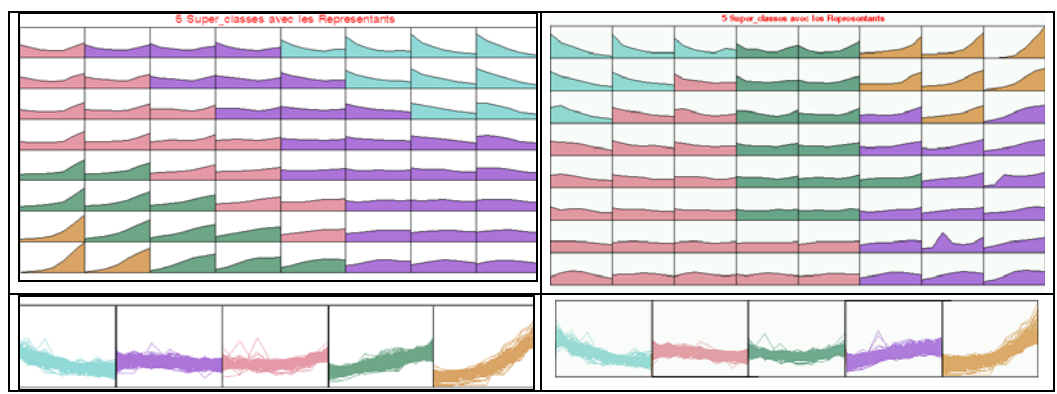

Fig. 6 : *Les communes, les classes et les super-classes,
initialisation II (dans les entrées),
à gauche, algorithme SOM, à droite, algorithme KBATCH*

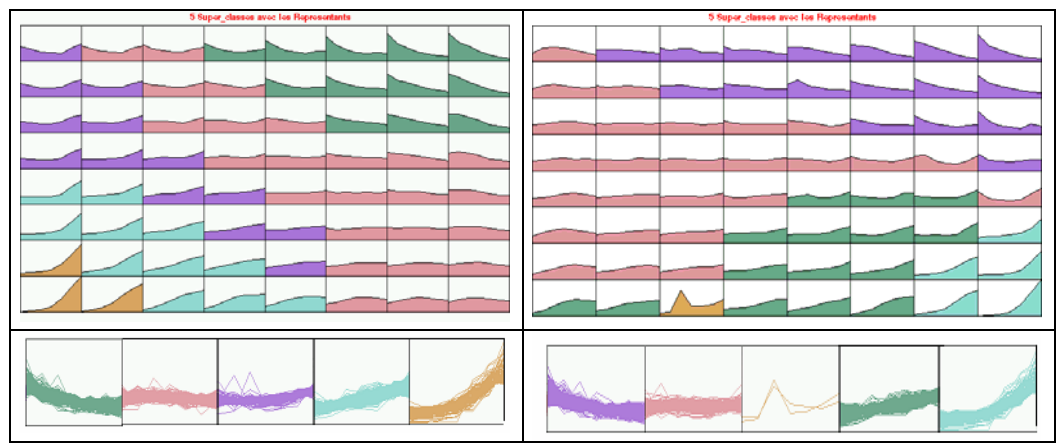

Fig. 7 : *Les communes, les classes et les super-classes,
initialisation III (dans le premier plan factoriel),
à gauche, algorithme SOM, à droite, algorithme KBATCH*

On remarque que l'algorithme KBATCH donne sensiblement les mêmes résultats que l'algorithme SOM pour l'initialisation II, alors que pour l'initialisation III, il fournit des classes d'effectifs très déséquilibrés. On voit donc que l'algorithme KBATCH est extrêmement sensible à l'initialisation, alors que l'algorithme SOM est beaucoup plus robuste. Cependant, dans la majorité des cas, on constate que l'algorithme KBATCH conduit à des minima plus petits de la fonction distorsion que l'algorithme SOM, au détriment de la qualité de l'organisation. D'autre part, l'algorithme SOM fournit une meilleure classification au sens des statistiques de Wilks (généralisation des statistiques de Fisher au cas multidimensionnel). Rappelons qu'une classification est d'autant meilleure que la valeur de la statistique de Wilks est faible. De même, les super-classes constituées à partir de l'algorithme SOM correspondent à des variances expliquées plus importantes, ce qui va dans le sens de classifications de meilleure qualité.

Le tableau ci-dessous indique pour chacun des six exemples les valeurs finales de la distorsion après calcul des vecteurs codes, les effectifs des super-classes, les statistiques de Wilks et les pourcentages d'inertie expliquée (rapport entre l'inertie inter-classes et l'inertie totale).

|  | Dist | Cl 1 | Cl 2 | Cl 3 | Cl 4 | Cl 5 | Wilks | % inert |
|---|---|---|---|---|---|---|---|---|
| **SOM I** | 0.0073 | 323 | 320 | 726 | 322 | 92 | **0.074** | **82.1%** |
| **BATCH I** | 0.0112 | 1783 | 0 | 0 | 0 | 0 |  | 43.1% |
| **SOM II** | 0.0072 | 322 | 700 | 351 | 315 | 95 | **0.078** | **82.2%** |
| **BATCH II** | **0.0059** | 213 | 653 | 304 | 402 | 211 | 0.089 | 73.8% |
| **SOM III** | 0.0073 | 402 | 676 | 286 | 325 | 94 | **0.083** | **81.8%** |
| **BATCH III** | **0.0058** | 588 | 640 | 2 | 404 | 149 | 0.102 | 71.9% |

## 8. Croisement avec des variables qualitatives

Pour interpréter les classes selon des variables qualitatives non utilisées dans l'algorithme de classement de Kohonen, il peut être intéressant d'étudier la répartition de leurs modalités dans chaque classe. Après avoir calculé des statistiques élémentaires dans chaque classe, on peut dessiner à l'intérieur de chaque cellule un camembert montrant comment sont réparties les modalités de chacune des variables qualitatives, comme le montrent les figures 8 et 9.

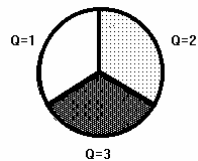

Fig. 8 : *Camembert, pour trois modalités équi-réparties.*

Dans la figure 9, il s'agit toujours des courbes de consommation (Exemple I) et la variable qualitative est le jour avec 3 niveaux : dimanche (noir), samedi (gris) et jours de semaine (blanc).

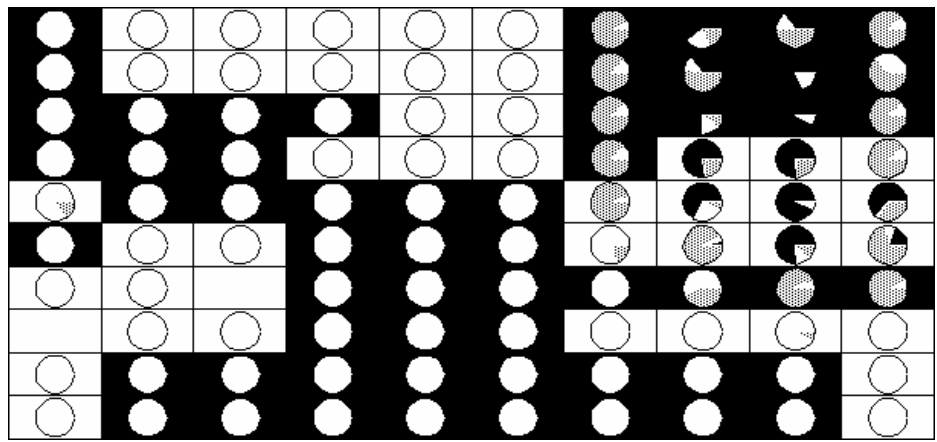

Fig. 9 : *Consommations électriques demi-horaires : dans chaque cellule, la répartition de la variable Jour est représentée. Les unités 8, 9, 18, 19, 28, 29, 38, 39, 48, 49, 50, 59 contiennent principalement des dimanches. Les unités 7, 17, 27, 37, 47, 58, 59, 10, 20, 30, 40, 60, 70, contiennent surtout des samedis. Les autres unités ne contiennent que des jours de semaines. Les cellules sont numérotées par ligne.*

Dans le cas de l'exemple des communes de la vallée du Rhône (Exemple II), on peut représenter sur la carte obtenue (figure 5, en haut à gauche) la répartition suivant la variable qualitative « Département ». Ces communes appartiennent à 8 départements différents (Ardèche, Bouches du Rhône, Drôme, Gard, Hérault, Isère, Haute-Loire, Vaucluse).

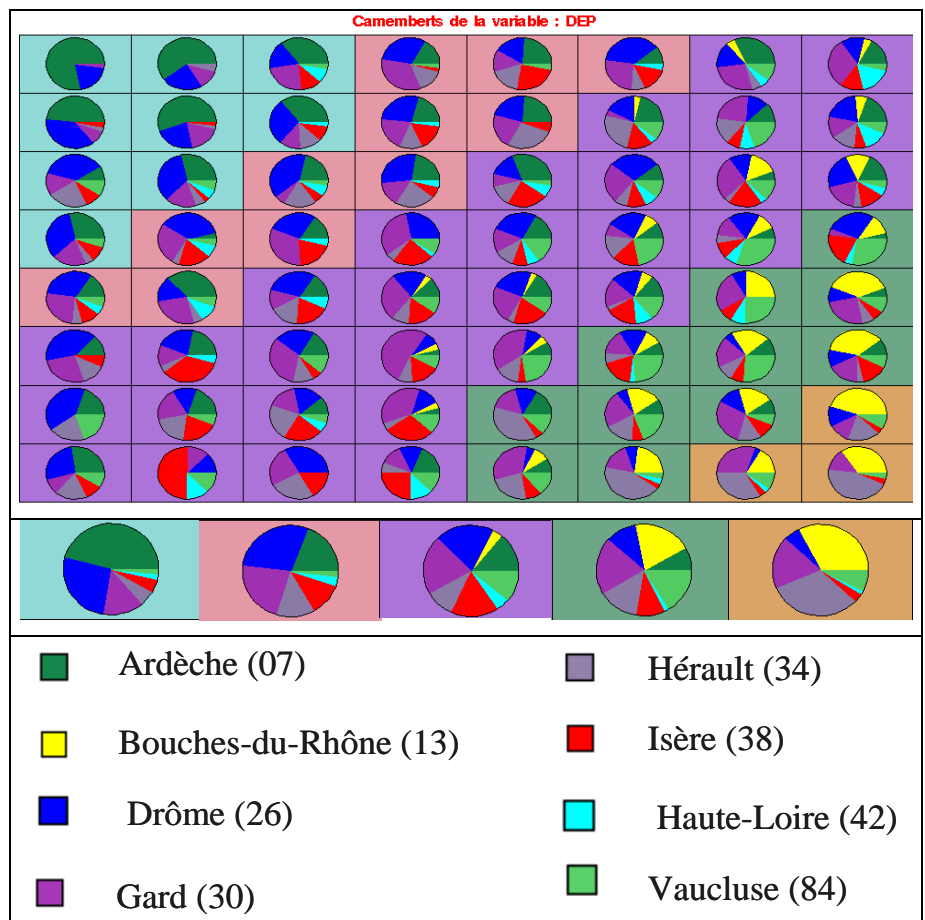

Fig. 10 : *Les communes : Répartition des 8 départements sur la carte de Kohonen de la figure 5 , et dans les super-classes correspondantes.*

On constate que les communes en forte décroissance appartiennent très majoritairement à l'Ardèche et à la Drôme, et ne comprennent aucune commune des Bouches du Rhône, alors que les communes des Bouches du Rhône sont presque toutes dans les super-classes 4 et 5, c'est-à-dire en croissance démographique.

On peut aussi interpréter la carte de Kohonen au moyen des variables quantitatives qui ont servi à la classification. On dessine alors le long de la carte ou de la ficelle les valeurs de chaque variable (figure 11).

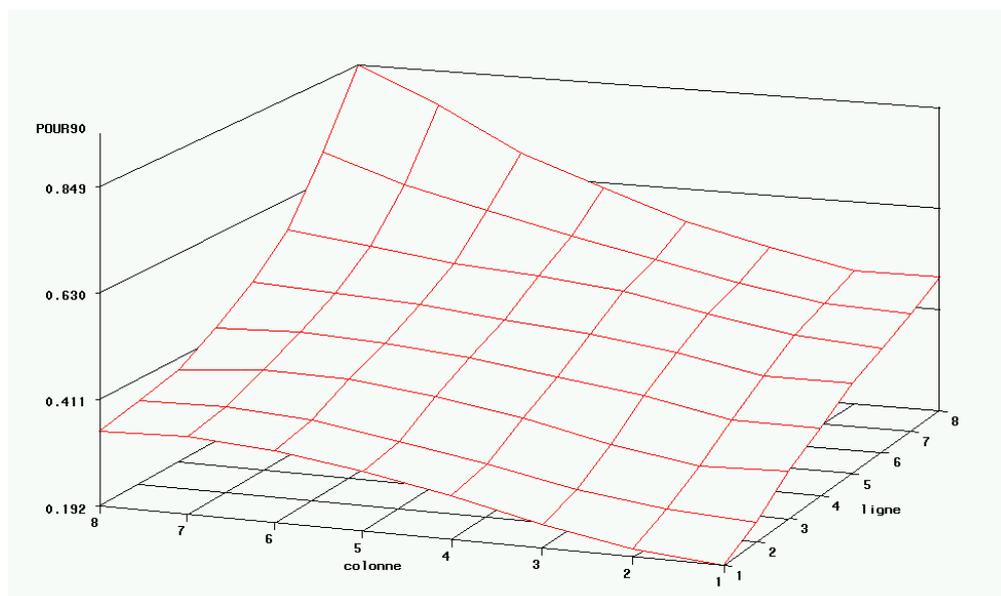

Fig. 11 : *Les communes : représentation de la variable « recensement en 1990 » le long de la carte de Kohonen de la figure 5.*

Ces cartes permettent de vérifier l'organisation, composante par composante, et de prolonger la descriptions des classes et super-classes.

On peut également les représenter de manière à mettre en évidence les super-classes (figure 12).

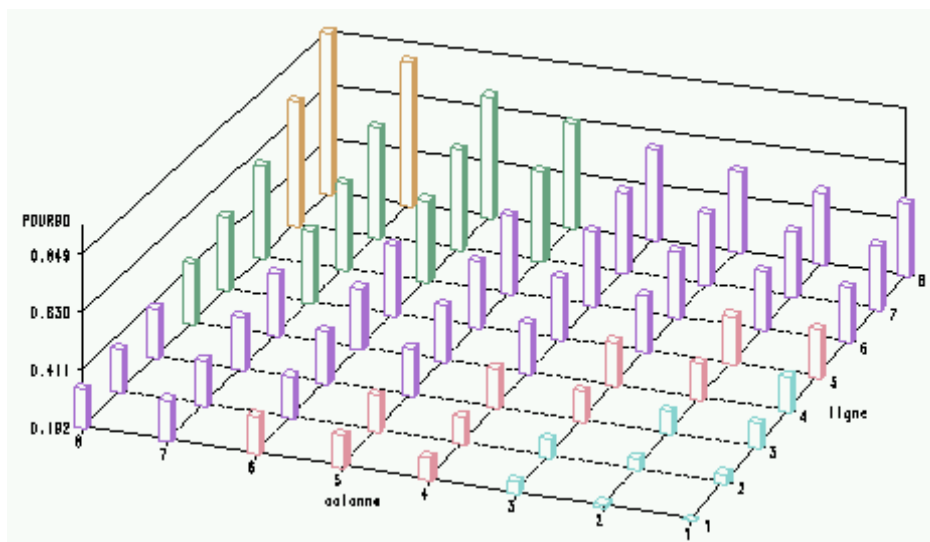

Fig. 12 : *Les communes : représentation de la variable « recensement en 1990 » le long de la carte de Kohonen de la figure 5, les couleurs correspondent aux super-classes.*

## 9. Analyse des relations entre variables qualitatives

On définit ici plusieurs algorithmes originaux qui permettent d'analyser les relations entre variables qualitatives. Voir Cottrell, Letrémy, Roy (1993), la thèse de Smaïl Ibbou (1998), et Cottrell, Letrémy (2001, 2003).

Attention, la plupart du temps les variables qualitatives ne peuvent pas être utilisées telles quelles, même lorsque les modalités sont codées par des nombres. S'il n'existe pas de relation d'ordre sur les codes (1 pour yeux bleus, 2 pour yeux marrons, etc.), cela n'a aucun sens de les utiliser comme des variables numériques pour faire un apprentissage de Kohonen. Même si les codes correspondent à une progression croissante ou décroissante, cela n'aurait un sens que si une échelle linéaire était utilisée (la modalité 2 correspondant à la moitié de la progression entre la modalité 1 et 3).

Si on souhaite simplement obtenir un classement des individus et non les relations entre modalités, une méthode utile consiste alors à effectuer préalablement au traitement des variables qualitatives une analyse des correspondances multiples classique, en gardant toutes les coordonnées, ce qui revient à coder tous les individus par les coordonnées résultant de cette transformation. Une fois les individus représentés par des variables numériques, on peut les classer par l'algorithme de Kohonen, mais on a perdu les modalités et d'autre part les calculs sont lourds et coûteux en temps, ce qu'on cherche justement à éviter en utilisant l'algorithme de Kohonen.

Nous allons dans la suite présenter plusieurs algorithmes et indiquer dans quels cas ils sont utiles. En effet, on peut traiter les données de deux façons : ou bien on s'intéresse seulement aux relations entre les modalités des variables qualitatives, ou bien on cherche à classer simultanément les modalités et les individus. Les méthodes proposées sont légèrement différentes dans l'un ou l'autre cas.

### 9.1 Les données et les notations

**Tableau disjonctif complet**

Pour tous les algorithmes, on considère $N$ individus et un certain nombre $K$ de variables qualitatives. La variable $k = 1,2, ..., K$ a $m_k$ modalités. Chaque individu choisit une et une seule modalité pour chaque variable. Si $M$ est le nombre total de modalités, chaque individu est représenté par un $M$-vecteur composé de 0 et de 1. Il n'y a qu'un 1 parmi les $m_1$ premières composantes, seulement un 1 entre la $m_1+1$-ième et la $(m_1+m_2)$-ième, etc. La matrice de données correspondante est appelé *tableau disjonctif complet*, on le notera $D$.

Si l'on ne s'intéresse pas aux réponses individuelles (ou si on n'en dispose pas) on résume les données en croisant deux à deux les K variables.

**Tableau de contingence**

S'il n'y a que deux variables ($K = 2$), on considère alors le tableau de contingence $T$ qui croise les deux variables. Dans ce cas, on note $p$ (resp. $q$) à la place de $m_1$ (resp. $m_2$). Dans la table de contingence $T$, la première variable qualitative a $p$ niveaux et correspond aux lignes. La seconde a $q$ niveaux et correspond aux colonnes. Le terme $n_{ij}$ est le nombre d'individus appartenant à la fois à la classe $i$ et à la classe $j$. A partir de la table de contingence, on calcule la matrice des fréquences relatives ($f_{ij} = n_{ij}/(\Sigma_{ij}\, n_{ij})$).

**Table de Burt**

Dans le cas général, quand $K > 2$, les données sont résumées dans une table de Burt qui est un tableau de contingence généralisé. C'est une matrice symétrique de dimension $M \times M$ composée de $K \times K$ blocs, tels que le bloc $B_{kl}$ (pour $1 \leq k \leq l \leq M$) est la table de contingence ($m_k \times m_l$) des variables $k$ et $l$. Le bloc $B_{kk}$ est une matrice diagonale, dont la diagonale est formée des nombres d'individus qui ont choisi les modalités 1, 2, ..., $m_k$, pour la variable $k$. Dans la suite, la table de Burt est notée $B$. On remarque que $B = D'D$.

## 9.2 Traitement des données

On peut traiter ces données de deux façons : ou bien on ne s'intéresse qu'aux relations entre les modalités des variables qualitatives, ou bien on cherche à classer simultanément les modalités et les individus. Les méthodes proposées sont légèrement différentes dans l'un ou l'autre cas.

Lorsque l'on ne s'intéresse qu'aux relations entre modalités, le cas le plus simple est celui où l'on dispose seulement de deux variables qualitatives. Alors, les données sont résumées dans un **tableau de contingence** et l'analyse des correspondances classique consiste à faire une analyse en composantes principales pondérée, utilisant la distance du Chi-deux, simultanément sur les profils lignes et les profils colonnes. Dans ce cas précis, on introduit un algorithme inspiré de l'algorithme de Kohonen et appelé KORRESP. Il n'est défini que pour l'analyse des relations entre deux variables qualitatives.

S'il y a plus de deux variables qualitatives, et qu'on ne s'intéresse qu'aux modalités, on utilisera l'algorithme KACM.

Les autres algorithmes (KACM1 et KACM2, KDISJ) sont dédiés aux cas où on s'intéresse non seulement aux modalités, mais aussi aux individus.

## 10. L'algorithme KORRESP

L'algorithme KORRESP n'est défini que pour l'analyse des relations entre les modalités de deux variables qualitatives, c'est-à-dire lorsque $K = 2$. Il n'utilise que la table de contingence $T$. On ne dispose pas des réponses individuelles.

A partir du tableau de contingence, on calcule d'abord les profils lignes $r(i)$, $1 \leq i \leq p$ (le profil $r(i)$ est la distribution conditionnelle de la seconde variable quand la première variable vaut $i$) et les profils colonnes $c(j)$, $1 \leq j \leq q$ (le profil $c(j)$ est la distribution conditionnelle de la première variable quand la seconde vaut $j$).

On construit alors une nouvelle matrice de données $T^{corr}$ : à chaque profil ligne $r(i)$, on associe le profil de la colonne $c(j(i))$ le plus probable sachant $i$, et réciproquement on associe à chaque profil colonne $c(j)$ le profil ligne $r(i(j))$ le plus probable sachant $j$. La matrice de données $T^{corr}$ est la matrice $((p+q) \times (q+p))$ dont les $p$ premières lignes sont les vecteurs $(r(i),c(j(i)))$ et les $q$ dernières sont les vecteurs $(r(i(j)),c(j))$. On applique ensuite l'algorithme de Kohonen aux lignes de cette matrice $T^{corr}$.

Il faut noter qu'on tire au hasard les entrées alternativement parmi les $p$ premières lignes et parmi les $q$ dernières, mais que l'unité gagnante est calculée seulement sur les $q$ premières composantes dans le premier cas et sur les $p$ dernières dans le second cas, et selon la distance du $\chi^2$. Après convergence, chaque modalité des deux variables est classée dans une classe de Kohonen. Des modalités « proches » sont classées dans la même classe ou dans des classes voisines. Cet algorithme est une méthode très rapide et efficace d'analyse des relations entre deux variables qualitatives. Comme en analyse des correspondances classiques, on peut interpréter les proximités entre les modalités des deux variables en terme de corrélations.

On pourra voir par exemple Cottrell et Letrémy (1993, 1994) pour des applications sur données réelles.

Les graphes ci-dessous correspondent à un tableau de contingence très simple, qui croise la variable « type de monuments » et la variable « propriétaire » (Exemple III).

La variable « type de monument » a 11 modalités : antiquités préhistoriques (preh), antiquités historiques (hist), châteaux (chat), architecture militaire (mili), cathédrales (cath), églises (egli), chapelles (chap), monastères (mona), édifices civils publics (ecpu), édifices civils privés (ecpr), divers (dive).

La variable « propriétaire » a 6 modalités : commune (COMM), privé (PRIV), état (ETAT), département (DEPA), établissement public (ETPU), non déterminé (NDET).

La table suivante est le tableau de contingence qui croise les deux variables ainsi définies.

**Monuments Historiques classés par catégorie et type de propriétaire**

| MONUMENT | COMMUNE | PRIVE | ETAT | DEPARTEMENT | ETAB PUBLICS | NON DETERMINE |
|---|---|---|---|---|---|---|
| préhistorique | 244 | 790 | 115 | 9 | 12 | 144 |
| historique | 246 | 166 | 46 | 23 | 11 | 31 |
| chateau | 289 | 964 | 82 | 58 | 40 | 2 |
| militaire | 351 | 76 | 59 | 7 | 2 | 0 |
| cathédrale | 0 | 0 | 87 | 0 | 0 | 0 |
| église | 4298 | 74 | 16 | 5 | 4 | 2 |
| chapelle | 481 | 119 | 13 | 7 | 8 | 4 |
| monastère | 243 | 233 | 44 | 37 | 18 | 0 |
| civil publique | 339 | 47 | 92 | 19 | 41 | 2 |
| civil privé | 224 | 909 | 46 | 7 | 18 | 4 |
| divers | 967 | 242 | 109 | 40 | 10 | 9 |

Les deux figures 13 et 14 montrent respectivement le premier plan factoriel obtenu par une analyse des correspondances classique (21% de variance expliquée), et une grille de Kohonen de taille 5 par 5 obtenue selon les méthodes décrites dans ce paragraphe.

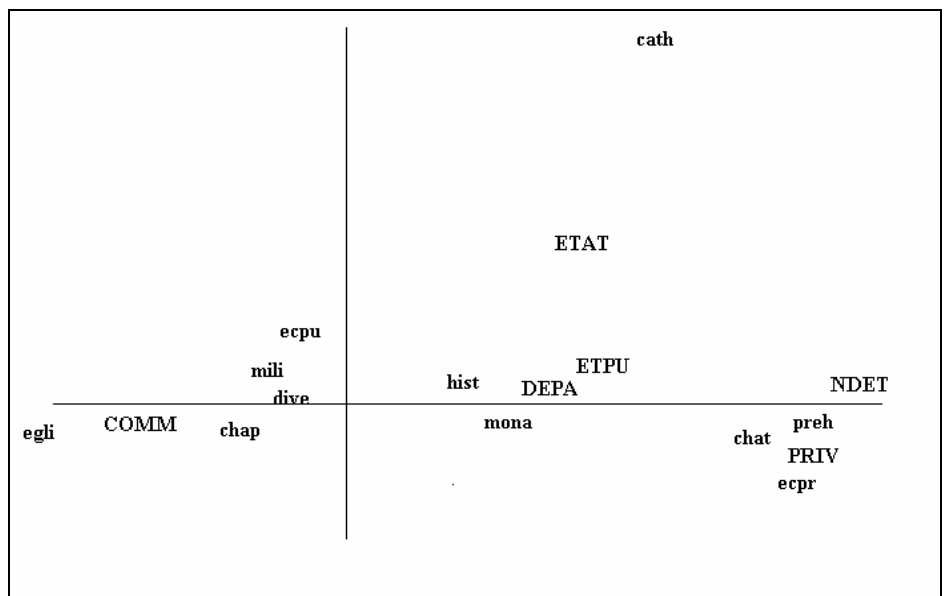

Fig. 13 : *Les monuments : premier plan factoriel de l'analyse des correspondances classiques.*

Fig. 14 : *Les monuments : classements des modalités des deux variables qualitatives. On observe que les cathédrales appartiennent à l'état, les châteaux sont privés, les églises et les chapelles sont à la commune, etc.*

Les associations sont bien les mêmes, mais elles sont plus faciles à voir sur une carte de Kohonen qui constitue une seule représentation. Il n'est pas nécessaire de regarder les projections successives comme dans une analyse classique.

## 11. L'algorithme KACM

Lorsque le nombre $M$ de modalités est supérieur (ou égal) à 2, et quand on ne cherche pas à classer les individus, mais seulement les modalités, on peut utiliser un autre algorithme inspiré de l'algorithme de Kohonen, appelé KACM.

On prend alors comme matrice des données la table de Burt $B$. Les lignes sont normalisées avec somme 1. A chaque étape, on tire une ligne normalisée au hasard selon la loi de probabilité donnée par la distribution empirique de la modalité correspondante. On définit l'unité gagnante selon la distance du $\chi^2$ et on met à jour les vecteurs codes comme d'habitude.

Cela est équivalent à utiliser un algorithme de Kohonen usuel (avec distance euclidienne et tirage uniforme des lignes) sur la matrice de Burt corrigée, où chaque terme $n_{ij}$ et remplacé par

$$n_{ij}^c = \frac{n_{ij}}{\sqrt{n_{i.}}\sqrt{n_{.j}}}$$

où $n_{i.}$ et $n_{.j}$ sont les sommes en lignes et en colonnes de la table de Burt. Remarquons que, puisque la matrice de Burt est symétrique, $n_{i.} = n_{.i}$ pour tout $i$.

Après convergence, on obtient une classification organisée des modalités, telle que des modalités « liées » appartiennent à la même classe ou à des classes voisines. Dans ce cas aussi, la méthode KACM fournit une technique alternative très intéressante à l'Analyse des Correspondances Multiples classique.

Pour illustrer cette méthode, prenons un exemple extrait de l'Enquête Emploi du temps de l'INSEE (1998-1999) : Exemple IV. On peut se reporter à Letrémy et coll. (2002), pour l'étude complète, Cottrell et Letrémy (2001) pour une première étude concernant les salariés à temps partiel. Les exemples présentés ici portent sur 207 salariés à temps partiel en CDD ou en CDI. Ils sont décrits par 8 variables qualitatives et 23 modalités définies dans la table suivante :

| Variables | Modalités | Noms |
|---|---|---|
| Type de contrat | CDI, CDD | CDI, CDD |
| Age | <25, [25, 40], [40,50], ≥50 | AGE1,AGE2,AGE3,AGE4 |
| Horaires de travail par jour | Identique, postés, variables | HOR1,HOR2,HOR3 |
| Travail le samedi | Jamais, parfois, habituellement | SAM1,SAM2,SAM3 |
| Travail le dimanche | Jamais, parfois, habituellement | DIM1,DIM2,DIM3 |
| Possibilité d'absence | Oui, oui sous condition, non | ABS1,ABS2,ABS3 |
| Temps partiel choisi | Oui, non | CHO1,CHO2 |
| Possibilité de récupération des heures | Sans objet, oui, non | REC0,REC1,REC2 |

Un simple tri montre que les contrats de type CDI représentent 83.57 % de la population, alors que les temps partiels subis (non choisis) représentent 46% des individus. On se servira de ces données pour illustrer les méthodes adaptées aux variables qualitatives (avec ou sans représentation simultanée des individus).

La figure 15 montre le résultat de l'analyse des correspondances multiples de ces données.

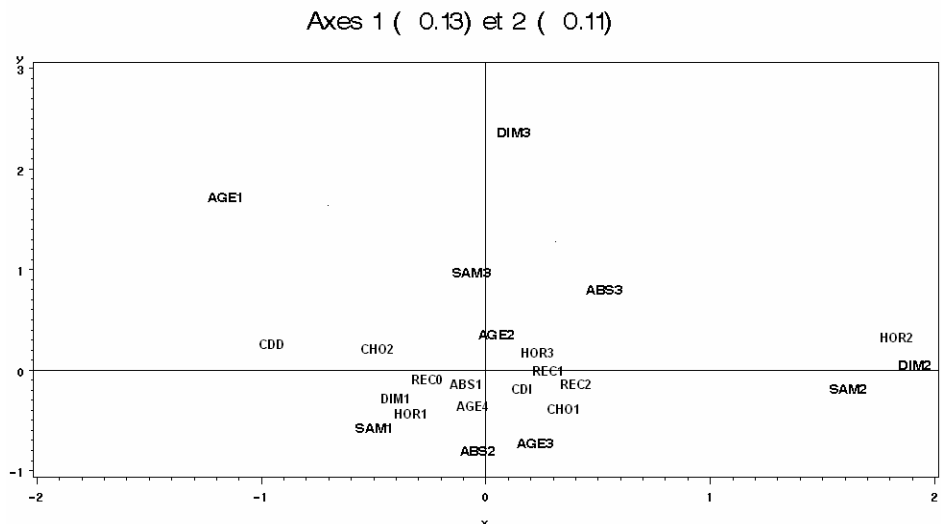

Fig. 15 : *Premier plan factoriel de l'analyse des correspondance multiples.*

Le résultat de l'algorithme KACM permet une représentation unique (figure 16) des modalités regroupées en classes.

| ABS2 REC2 | | SAM2 DIM2 | ABS3 | HOR2 DIM3 |
|---|---|---|---|---|
| | | | | |
| AGE3 | | AGE4 REC1 | | CDD AGE1 |
| CHO1 | SAM1 | | CHO2 | |
| CDI HOR1 | | DIM1 ABS1 | REC0 | AGE2 HOR2 SAM3 |

Fig. 16 : *Carte de Kohonen des modalités (KACM).*

On remarque les mêmes regroupements que sur la figure précédente.

On peut adapter de deux façons cet algorithme au cas où l'on cherche à classer également les individus. Nous suivons ici les propositions de Smaïl Ibbou dans sa thèse (Ibbou, 1998).

**Algorithme KACM2**

On peut procéder en pratiquant comme dans KACM un algorithme de Kohonen sur la matrice de Burt corrigée (voir plus haut). Puis on classe les individus comme des données supplémentaires, après les avoir normalisés, en remplaçant

$$d_{ij} \text{ par } \frac{d_{ij}}{K}, \text{ pour } j = 1, \cdots, M.$$

Cette transformation permet de les plonger dans le même espace que les lignes de la matrice de Burt corrigée $B^c$.

Avec KACM2, l'algorithme est rapide puisqu'il ne porte que sur les *M* modalités, mais il faut en général plus d'itérations que pour réaliser simplement une classification des modalités au moyen de KACM. Cela est dû au fait qu'il faut non seulement bien classer les modalités, mais aussi déterminer avec précision les vecteurs codes qui servent ensuite à classer les individus en données supplémentaires.

**Algorithme KACM1**

On peut procéder « à l'envers ». Puisqu'on cherche à représenter simultanément les modalités et les individus, on considère alors le tableau disjonctif complet *D*. On le corrige de manière convenable (pour introduire la distance du Chi-deux, en posant

$$d_{ij}^c = \frac{d_{ij}}{\sqrt{d_{i.}}\sqrt{d_{.j}}},$$

où $d_{i.} = K$ est la somme en ligne et $d_{.j}$ est la somme en colonne (effectif de la modalité *j*).

On pratique un algorithme de Kohonen classique sur ce tableau de données corrigées, on classe ainsi les *N* individus. On classe ensuite les modalités normalisées comme des « individus-types » en données supplémentaires.

La modalité *j* est représentée par le *M*-vecteur de coordonnées :

$$\frac{n_{jl}}{d_{.j}\sqrt{d_{.l}}\sqrt{K}}, \text{ pour } l = 1,\cdots, M.$$

La représentation graphique n'est pas toujours utile (quand il y a trop d'individus), mais la classification permet d'associer à chaque modalité, ou groupe de modalités, un ensemble d'individus particulièrement représentatifs de cette modalité ou groupe de modalités. Voir Cottrell et Ibbou (1995), Cottrell, de Bodt et Henrion (1996) pour plus de détails et pour des applications à des données financières.

Les deux méthodes précédentes permettent la représentation simultanée des individus et des modalités, mais rompent la symétrie. Au contraire, l'algorithme introduit ci-dessous (noté KDISJ) conserve cette symétrie et est directement inspiré des méthodes classiques. Il a été introduit dans Cottrell et Letrémy (2001, 2003).

## 12. L'algorithme KDISJ

Il s'agit de classer *simultanément* les individus et les modalités des variables qualitatives qui les décrivent.

Avec les notations introduites, on considère ici le tableau disjonctif complet, noté *D*. Remarquons qu'il contient toute l'information permettant de connaître aussi bien les individus que la répartition des modalités. On note $d_{ij}$ le terme général de ce tableau qui peut être considéré comme un tableau de contingence croisant la variable "individu" à *N* modalités et la variable "modalités" à *M* modalités. Le terme $d_{ij}$ prend ses valeurs dans {0,1}.

On utilise alors une adaptation de l'algorithme KORRESP présenté plus haut et qui porte maintenant sur ce « tableau de contingence ».

De manière à utiliser la distance du $\chi^2$ sur les lignes autant que sur les colonnes, et pour pondérer les modalités de façon proportionnelle à leurs effectifs, on corrige le tableau disjonctif complet, et on pose

$$d_{ij}^c = \frac{d_{ij}}{\sqrt{d_{i.}d_{.j}}}$$

où :

$$d_{i.} = \sum_{j=1}^{M} d_{ij} \text{ et } d_{.j} = \sum_{i=1}^{N} d_{ij}.$$

Remarquons que, dans le cas d'un tableau disjonctif complet, $d_{i.}$ vaut *K*, quelque soit *i*. Le terme $d_{.j}$ est l'effectif de la modalité *j*.

Le tableau ainsi corrigé est noté $D^c$ (tableau disjonctif corrigé). Cette transformation est la même que celle qui a été proposée par Smaïl Ibbou dans sa thèse (Ibbou, 1998, Cottrell et Ibbou, 1995). Après cette transformation, utiliser la distance euclidienne sur $D^c$ équivaut à utiliser la distance du $\chi^2$ pondérée sur *D*.

Ces corrections sont exactement celles qu'on fait traditionnellement lorsqu'on pratique une analyse des correspondances, qui équivaut en fait à une analyse en composantes principales pondérée, utilisant la distance du Chi-deux, simultanée sur les profils lignes et les profils colonnes.

On choisit ensuite un réseau de Kohonen, et on associe à chaque unité *u* un vecteur code $C_u$ formé de (*M* + *N*) composantes, les *M* premières évoluent dans l'espace des individus (représentés par les lignes de $D^c$), les *N* dernières dans l'espace des modalités (représentées par les colonnes de $D^c$). On note

$$C_u = (C_M, C_N)_u = (C_{M,u}, C_{N,u})$$

pour mettre en évidence la structure du vecteur code $C_u$. Les étapes de l'apprentissage du réseau de Kohonen sont doubles. On tire alternativement une ligne de $D^c$ (c'est-à-dire un individu *i* ), puis une colonne (c'est-à-dire une modalité *j* ).

Quand on tire un individu *i*, on lui associe la modalité *j(i)* définie par

$$j(i) = Arg \max_{j} d_{ij}^{c} = Arg \max_{j} \frac{d_{ij}}{\sqrt{Kd_{.j}}}$$

qui maximise le coefficient $d_{ij}^{c}$, c'est-à-dire la modalité la plus rare dans la population totale parmi les modalités qui lui correspondent. Ensuite on crée un vecteur individu étendu $X = (i, j(i)) = (X_M, X_N)$, de dimension $(M + N)$ (voir figure 17). On cherche alors parmi les vecteurs-codes celui qui est le plus proche, au sens de la distance euclidienne restreinte aux $M$ premières composantes. Notons $u_0$ l'unité gagnante. On rapproche alors les vecteurs-codes de l'unité $u_0$ et de ses voisines du vecteur complété $(i, j(i))$, selon la loi usuelle de Kohonen.

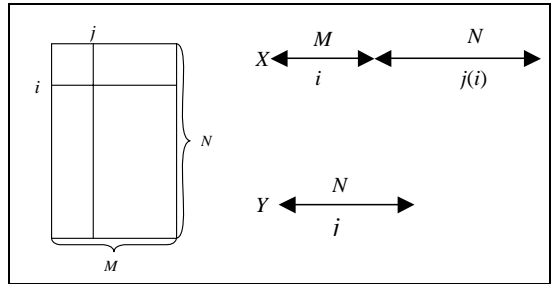

*Fig.* 17 : La matrice $D^c$, les vecteurs $X$ et $Y$.

On peut écrire formellement cette étape ainsi :
$$\begin{cases} u_0 = Arg \min_{u} \|X_M - C_{M,u}\| \\ C_u^{new} = C_u^{old} + \varepsilon\, \sigma(u, u_0)(X - C_u^{old}) \end{cases}$$

où ε est le paramètre d'adaptation (positif, décroissant), et σ est la fonction de voisinage telle que $\sigma(u, u_0)$ vaut 1 si $u$ et $u_0$ sont voisines sur le réseau de Kohonen, et 0 sinon.

Quand on tire une modalité $j$, de dimension $N$ (une colonne de $D^c$), on ne lui associe pas de vecteur, en effet, par construction, il y a beaucoup d'ex-æquo et le choix serait arbitraire. On cherche parmi les vecteurs codes celui qui est le plus proche, au sens de la distance euclidienne restreinte aux $N$ dernières composantes. Soit $v_0$ l'unité gagnante. On rapproche alors les $N$ dernières composantes du vecteur-code gagnant associé à $v_0$ et de ses voisins de celles du vecteur modalité $j$, sans modifier les $M$ premières composantes. Pour simplifier notons $Y$ (voir figure 17) le vecteur colonne de dimension $N$ correspondant à la modalité $j$. Cette étape peut s'écrire :
$$\begin{cases} v_0 = Arg \min_{u} \|Y - C_{N,u}\| \\ C_{N,u}^{new} = C_{N,u}^{old} + \varepsilon\, \sigma(u, u_0)(Y - C_{N,u}^{old}) \end{cases}$$

où les $M$ premières composantes ne sont pas modifiées.

On pratique ainsi un classement classique de Kohonen sur les individus, un autre sur les modalités, tout en les maintenant associés. Après convergence, les individus et les modalités sont classés dans les classes de Kohonen. Des individus ou modalités « proches » sont classés dans la même classe ou dans des classes voisines. On appelle KDISJ l'algorithme ainsi défini. Il suffit en général de faire $15(M+N)$ itérations pour obtenir la convergence.

## 13. Comparaison des cartes ainsi obtenues

Nous pouvons comparer les cartes obtenues à partir de l'exemple IV (Enquête Emploi du temps).

On donne figure 18 le premier plan factoriel de l'analyse des correspondances où sont représentés simultanément les modalités et les individus.

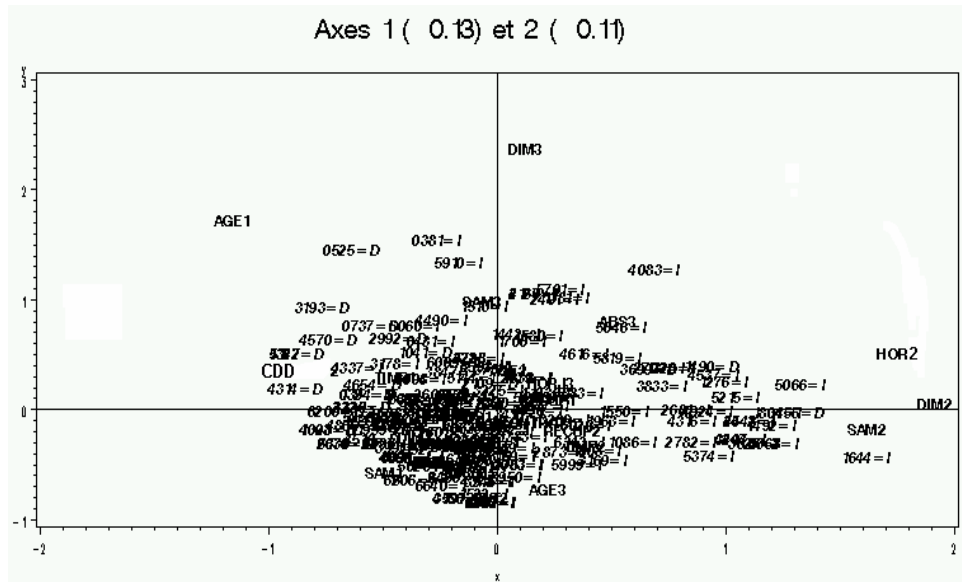

Fig. 18 : *Premier plan factoriel avec les individus.*

On constate que ce genre de carte est plutôt illisible et qu'il est pratiquement impossible d'en déduire une classification.

On représente ensuite les cartes de Kohonen obtenues par les méthodes KACM2 (figure 19), KACM1 (figure 20) et KDISJ (figure 21). Entre parenthèses, nous indiquons le nombre d'individus en CDI, et le nombre d'individus en CDD.

| ABS2 REC2 |       | SAM2 DIM2 | ABS3   | HOR2 DIM3 |
|-----------|-------|-----------|--------|-----------|
| (1,0)     | (1,0) | (14,3)    | (2,0)  | (1,0)     |
|      |      | (3,0)        |      |             |
| AGE3 |      | AGE4<br>REC1 |      | CDD<br>AGE1 |
| (21,0) | (1,0) | (6,2)      |      | (0,3)       |
| CHO1 | SAM1 |              | CHO2 |             |
| (6,0) | (47,2) | (1,0)     | (13,22) | (1,2)    |
| CDI<br>HOR1 |  | DIM1<br>ABS1 | REC0 | AGE2<br>HOR2<br>SAM3 |
| (7,0) | (1,0) | (9,0)      | (12,0) | (26,0)    |

Fig. 19: *Carte de Kohonen des modalités et des individus (KACM2).*

|      |      |      |      |      |
|------|------|------|------|------|
| HOR2<br>SAM2<br>DIM2 |  | REC2 |  | ABS2 |
| (11,2) | (11,0) | (3,0) | (4,0) | (16,1) |
|      | CDI, AGE3<br>HOR3<br>CHO1<br>REC1 |  |  |  |
| (8,1) | (8,0) | (16,1) | (2,0) | (3,0) |
|      |      | HOR1<br>ABS1 |      | AGE2 |
| (12,1) | (9,0) | (3,0) |      | (24,0) |
| SAM3 | AGE1<br>CHO2 |  | AGE4 | SAM1<br>DIM1<br>REC0 |
|      | (5,0) | (4,0) | (15,0) | (0,1) |
| DIM3 | ABS3 |      |      | CDD  |
| (12,1) |     | (7,1) | (0,6) | (0,19) |

Fig. 20: *Carte de Kohonen des modalités et des individus (KACM1).*

|       |       |       |       |       |
|-------|-------|-------|-------|-------|
| ABS2  |       | AGE2  | CHO2  | CDD   |
| (14,1) | (5,0) | (20,0) | (1,1) | (0,21) |

|        | AGE3   | CDI<br>HOR1<br>SAM1<br>CHO1 | HOR3<br>DIM1<br>ABS1<br>REC0 |        |
|--------|--------|------|------|--------|
| (1,0)  | (18,0) |      |      |        |
| REC1   |        | AGE4 |      | AGE1   |
| (14,0) | (3,0)  | (19,0) |    | (4,7)  |
| SAM2   | REC2   |      | SAM3 |        |
| (7,0)  | (3,0)  | (2,0) | (17,0) | (2,0) |
| HOR2   | DIM2   | ABS3 |      | DIM3   |
| (13,2) | (9,1)  | (11,1) |    | (10,0) |

Fig. 21: *Carte de Kohonen des modalités et des individus (KDISJ).*

On constate que les associations principales se retrouvent dans les trois cartes. Les modalités correspondant aux conditions de travail « normales » (HOR1, SAM1, DIM1, ABS1, CHO1) sont proches, et proches des AGE3 et AGE4, etc. On a colorié en grisé les classes où les CDD sont majoritaires. On voit que sur les cartes KACM2 et KDISJ, elles correspondent bien aux modalités CDD et AGE1 (les jeunes), alors que sur la carte KACM1, l'association est moins claire.

Plus précisément, sur ces cartes qui comprennent aussi les individus, il est possible de contrôler si les modalités sont correctement placées par rapport aux individus, en calculant leurs déviations. Rappelons que la déviation pour une modalité *m* (partagée par $n_m$ individus) et pour une classe *k* (avec $n_k$, individus) peut être définie comme la différence entre le nombre d'individus qui ont cette modalité *m* et appartiennent à la classe *k* et le nombre *théorique* $n_m n_k / n$ qui correspondrait à une répartition conforme à la distribution dans la population totale. Sur les cartes obtenues par KACM2 et KDISJ, on observe (en les calculant) que toutes les déviations sont positives, ce qui exprime que les modalités se trouvent dans des classes où elles sont majoritairement représentées. Pour la carte KACM1 (figure 20) une seule modalité a une déviation négative, il s'agit de HOR3 (horaires variables), qui serait mieux placée dans la classe voisine (en dessous et à gauche).

On peut ensuite procéder à des regroupements en super-classes, en utilisant une classification hiérarchique sur les vecteurs codes, et on peut là aussi contrôler les représentations obtenues en calculant les déviations à l'intérieur des super-classes.

Les principaux avantages des méthodes que nous venons de présenter sont leur simplicité, leur rapidité et leur faible temps calcul. Elles produisent une seule carte alors que les analyses classiques fournissent plusieurs représentations d'information décroissante. Elles sont plus grossières, mais permettent une interprétation rapide.

## 14. Conclusion et perspectives

Nous proposons un ensemble de méthodes pour analyser des données multidimensionnelles, en complément des méthodes linéaires classiques, quand les observations sont décrites par des variables quantitatives et qualitatives. Nous avons réalisé de nombreuses études à l'aide de ces méthodes : voir la bibliographie.

En fait, dans la pratique, il est nécessaire de combiner les différentes techniques entre elles.

Dans le cas de données quantitatives, on peut au préalable réduire la dimension en faisant une analyse en composantes principales et en ne conservant qu'un nombre réduit de coordonnées.

On peut se servir de la classification en super-classes réalisée à partir de l'algorithme de Kohonen pour définir une nouvelle variable qualitative et pratiquer une Analyse des Correspondances Multiples ou un KACM sur l'ensemble des variables qualitatives (les variables d'origine et la variable de classe qu'on vient de définir). Cela permet d'obtenir une typologie des classes et d'aider à l'interprétation.

On peut transformer les données qualitatives par une analyse factorielle des correspondances où l'on conserve tous (ou presque tous) les facteurs, puis classer les observations par un algorithme de Kohonen.

On peut utiliser une ficelle de Kohonen soit directement sur les données de départ, soit sur les vecteurs codes, lorsqu'on désire définir un score.

On peut également utiliser la classification de Kohonen combinée avec la prédiction de paramètres de niveau et variance pour prévoir des courbes (de consommation par exemple, voir Cottrell, Girard, Rousset, 1998).

## Références bibliographiques

**Annexe A :**

**Encadré 1 : Les choix à faire pour la mise en oeuvre de l'algorithme de Kohonen**

- Nombre de classes, de super-classes
- Géométrie de la carte : grille, ficelle, cylindre, tore, maillage hexagonal
- Initialisation des vecteurs codes (parmi les données, dans l'enveloppe convexe, dans le premier plan principal)
- Pour constituer les super-classes, on peut faire une classification hiérarchique des vecteurs codes, ou utiliser une ficelle de Kohonen sur les vecteurs codes (Cela permet une typologie plus facile à décrire, puisqu'une ficelle fournit un score ordonné.
- Le choix d'une ficelle sur les données de départ fournit un *score*, croissant ou décroissant.

**Encadré 2 : Les sorties utiles**

- Classification des données
- Représentations des classes de Kohonen, leurs contenus, les distances mutuelles
- Les vecteurs-codes
- Les super-classes, leurs contenus
- Statistiques mono et multidimensionnelles permettant de qualifier les classifications obtenues
- Visualisation de la distorsion étendue
- Variations de chaque variable selon les classes obtenues
- Fabrication d'une variable qualitative (numéro de la classe ou de la super-classe
- Croisement possible avec les variables qualitatives

**Annexe B**

**Algorithmes pour variables quantitatives, leurs noms dans les programmes téléchargeables sur le site du SAMOS**

| *Algorithmes* | *Procédures* | *Déterministe* | *Stochastique* | *Voisinage* | *Organisation* | *Initialisation dans le domaine* | *Initialisation dans les entrées* | *Initialisation dans le plan factoriel* |
|---|---|---|---|---|---|---|---|---|
| FORGY | **FASTCLUS** | * | | Non | Non | * | * | * |
| SCL | **KFAST** | | * | Non | Non | * | * | * |
| SOM | **KACP** | | * | Oui | Oui | * | * | * |
| Batch | **KBATCH** | * | | Oui | Oui | | * | * |

**Annexe C**

**Méthodes pour variables qualitatives**

1) On peut croiser les variables qualitatives avec les classes obtenues

2) Pour un *tableau de contingence* : on utilise **KORRESP**

3) Pour les seules modalités d'un tableau de réponses de plus de 2 questions (*table de Burt*), on utilise **KACM**

4) Pour les modalités et les individus d'un tableau de réponses de plus de 2 questions (*tableau disjonctif complet*) : **KACM2, KACM1, KDISJ**

• **KACM2** classe les modalités comme **KACM,** puis classe les individus comme des *modalités supplémentaires*, à partir du tableau disjonctif complet.

• **KACM1** classe d'abord les individus à partir du tableau disjonctif complet, puis les modalités comme des *individus supplémentaires*, calculés à partir de la table de Burt.

• **KDISJ** classe simultanément les individus et les modalités, à partir du tableau disjonctif complet.

**Annexe D**

**Données manquantes**

Deux possibilités
- Pendant l'apprentissage,

On se sert des observations avec données incomplètes comme des autres (tirées aléatoirement) et, pour calculer le vecteur code gagnant, on calcule les distances en se restreignant aux composantes présentes.
- Après l'apprentissage

On fait l'apprentissage avec les observations complètes, puis on classe les observations incomplètes dans les classes obtenues. Les distances sont calculées avec les composantes présentes.

**Annexe E**

**Références des quatre exemples utilisés**

**Exemple I, Consommation électrique**

Cottrell, M., Girard, B., Girard, Y., Muller, C. & Rousset, P., (1995) : Daily Electrical Power Curves : Classification and Forecasting Using a Kohonen Map, *From Natural to Artificial Neural Computation, Proc. IWANN'95*, Springer, p. 1107-1113.

Cottrell M., Girard B. et Rousset P. (1998) : Forecasting of curves using a Kohonen classification, *Journal of Forecasting*, 17, p. 429-439.

Rousset P. (1999) : Applications des algorithmes d'auto-organisation à la classification et à la prévision, Thèse Université Paris 1.
\*\*\*\*\*\*\*\*\*\*\*\*\*\*\*\*\*\*\*\*\*\*\*\*\*\*\*\*\*\*\*\*\*\*\*\*\*\*\*\*\*\*\*\*\*\*\*\*\*\*\*\*\*\*\*\*\*\*\*\*\*\*\*\*

**Exemple II, Les communes de la vallée du Rhône**

Cottrell, M., Letrémy, P. (1994) : Classification et analyse des correspondances au moyen de l'algorithme de Kohonen : application à l'étude de données socio-économiques, *Proc. Neuro-Nîmes*, 74-83.

Cottrell M., Gaubert P., Letrémy P., Rousset P., : Analyzing and representing multidimensional quantitative and qualitative data: Demographic study of the Rhone valley. The domestic consumption of the Canadian families, *WSOM'99,* In: Oja E., Kaski S. (Eds), *Kohonen Maps*, Elsevier, Amsterdam, 1-14, 1999.
\*\*\*\*\*\*\*\*\*\*\*\*\*\*\*\*\*\*\*\*\*\*\*\*\*\*\*\*\*\*\*\*\*\*\*\*\*\*\*\*\*\*\*\*\*\*\*\*\*\*\*\*\*\*\*\*\*\*\*\*\*\*\*\*

**Exemple III, Les monuments**

Cottrell, M., Letrémy, P. & Roy, E. (1993) : Analyzing a contingency table with Kohonen maps : a Factorial Correspondence Analysis, *Proc. IWANN'93*, J.Cabestany, J. Mary, A. Prieto (Eds.), Lecture Notes in Computer Science, Springer-Verlag, 305-311.

Cottrell, M. & Letrémy, P. (1994) : Classification et analyse des correspondances au moyen de l'algorithme de Kohonen : application à l'étude de données socio-économiques, *Proc. Neuro-Nîmes*, 74-83.
\*\*\*\*\*\*\*\*\*\*\*\*\*\*\*\*\*\*\*\*\*\*\*\*\*\*\*\*\*\*\*\*\*\*\*\*\*\*\*\*\*\*\*\*\*\*\*\*\*\*\*\*\*\*\*\*\*\*\*\*\*\*\*\*

**Exemple IV, Les travailleurs à temps partiel**

Cottrell M., Letrémy P., (2001, 2003) : Working times in atypical forms of employment : the special case of part-time work, *Conf. ACSEG 2001*, Rennes, 2001, to appear in *Connectionist Approaches in Economics and Management Sciences*, Lesage C. et Cottrell M., Eds, Kluwer, 2003.

Letrémy P., Cottrell M., Macaire S., Meilland C., Michon F., (2001, 2002) : Le temps de travail des formes particulières d'emploi, Rapport final, IRES, Noisy-le-Grand, February 2001, *Economie et Statistique*, Octobre 2002.